\documentclass [a4paper,twoside,12pt]{article}

\usepackage[utf8]{inputenc}

\usepackage{amsmath,amsfonts,amssymb}
\usepackage{vmargin,graphicx,theorem}
\usepackage[english]{babel}
\usepackage{enumerate}
\usepackage{color}
\usepackage{pst-fill,pst-grad,pst-plot,pst-eucl,pstricks-add,pst-node}

\setpapersize[portrait]{A4}
\setmarginsrb{1.5cm}{1cm}{1.5cm}{2.5cm}{1.5cm}{0cm}{0.5cm}{2cm}

\newcommand{\disp}{\displaystyle}
\newcommand{\rr}{\ensuremath{\mathbb{R}}}

\newcommand{\dR}{\ensuremath{\mathbb{R}}}

\newtheorem{ethm}{Theorem}[section]

\newtheorem{ecor}[ethm]{Corollary}

\newtheorem{eprop}[ethm]{Proposition}

\newtheorem{elem}[ethm]{Lemma}

\newtheorem{edefi}[ethm]{Definition}

\newtheorem{erem}[ethm]{Remark}

\newcommand{\proofend}{~$\rhd$}
\newcommand{\proofbegin}{~$\lhd$}

\newenvironment{eproof}
               {\noindent {\emph{\textbf{Proof}}}\\\proofbegin~}
               {\proofend\\}

\newcommand{\p}[4]{{#3}\!\left#1{#4}\right#2}

\newcommand{\PAR}[1]{\ensuremath{{\left(#1\right)}}} 

\renewcommand{\phi}{\varphi}
\newcommand{\ep}{{\varepsilon}} 

\renewcommand{\geq}{\geqslant}

\newcommand{\entf}[1]{{\rm{Ent}}_{#1}}

\newcommand{\ent}[2]{\p(){\entf{#1}}{#2}}

\def\disp{\displaystyle}

\newcommand{\R}{\dR}

\newcommand{\1}{\hbox{1}\!\!\hbox{I}}

\newcommand{\beq}{\begin{equation}}\newcommand{\eeq}{\end{equation}}
\begin{document}

\title{Dimensional contraction via Markov transportation distance}
\author{ Fran{\c c}ois Bolley\thanks{Ceremade, Umr Cnrs 7534, Universit\'e Paris-Dauphine, Place de Lattre de Tassigny, F-75775 Paris cedex 16. bolley@ceremade.dauphine.fr}, Ivan Gentil\thanks{Institut Camille Jordan, Umr Cnrs 5208, Universit\'e Claude Bernard Lyon 1, 43 boulevard du 11 novembre 1918, F-69622 Villeurbanne cedex. gentil@math.univ-lyon1.fr}\, and Arnaud Guillin\thanks{Institut Universitaire de France and Laboratoire de Math\'ematiques, Umr Cnrs 6620, Universit\'e Blaise Pascal, Avenue des Landais, F-63177 Aubi\`ere cedex. guillin@math.univ-bpclermont.fr}}

\date{\today}

\maketitle

\abstract{It is now well known that curvature conditions {\it \`a la} Bakry-\'Emery are equivalent to contraction properties of the heat semigroup with respect to the classical quadratic Wasserstein distance. However, this curvature condition may include a dimensional correction which up to now had not induced any strenghtening of this contraction. We first consider the simplest example of the Euclidean heat semigroup, and prove that indeed it is so. To consider the case of a general Markov semigroup, we introduce a new distance between probability measures, based on the semigroup, and adapted to it.  We prove, in the setting of a compact Riemannian manifold, that this Markov transportation distance satisfies the same properties as the Wasserstein distance does in the specific case of the Euclidean heat semigroup, namely dimensional contraction properties and  Evolution Variational Inequalities.
}

\bigskip

\bigskip


\section{Introduction}

Contraction properties of (Markov) semigroups are an important probabilistic and analytic tool: for instance they enable to study the existence of invariant probability measures, or the stability and long time behaviour of solutions to various linear (Fokker-Planck, kinetic Fokker-Planck,~...) or non linear (McKean-Vlasov, porous medium, Boltzmann,~...) partial differential equations. An important aspect is of course the distance in which we measure this contraction.  Recent progress has shown that the Wasserstein distance  is a particularly relevant and natural choice, in particular, but not only, for dynamics which have been interpreted as gradient flows for this distance (see for example \cite{otto,CT03,cmcv-06,cgm-08,bgm,nps,BGG11} and the reference books \cite{ambrosio-gigli-savare,villani-book1}). Here and below the Wasserstein distance between two Borel probability measures $\nu$ and $\mu$ on a Polish metric space $(E,d)$   is defined by 
$$
W_2(\mu,\nu)=\inf\Big(\int d^2(x,y)d\pi(x,y)\Big)^{1/2},
$$
where the infimum runs over all probability measures $\pi$ on $E\times E$ with marginals $\mu$ and $\nu$. We refer again to~\cite{ambrosio-gigli-savare,villani-book1} for a reference presentation of this distance and its interplay with the optimal transportation problem.

\medskip

On the other hand, geometric properties of metric spaces are an important and vast topic with many diverse issues, and the Wasserstein distance has provided new insight on them, see \cite{otto05,sturm-vonrenesse,sturm,LV,AGS13,AGS12}. A particularly relevant notion is the one of curvature which has recently attracted much attention. It turns out that it can be handled in terms of a contraction property in Wasserstein distance, as follows.

 Let $(H_t)_{t \geq 0}$ denote the heat semigroup on a smooth and complete (and connex) Riemannian manifold $(M,g)$: it solves the heat equation $\partial_t u=\Delta_g u$ where $\Delta_g$ is the Laplace-Beltrami operator on $M$. Let also $\mu$ be the Riemannian measure on $(M,g)$ and $d$ the associated Riemannian distance.
Then a fundamental result, due to M.~von Renesse and  K.-T.~Sturm in~\cite{sturm-vonrenesse}, says that the Ricci curvature of the manifold is bounded from below by a constant $R\in\R$ if and only~if 
\begin{equation}
\label{ine-vonrenesse}
W_2(H_t f\mu,H_tg\mu)\leq e^{-R t} W_2( f\mu,g\mu)
\end{equation}
for any $t\geq0$ and any probability densities $f,g$ with respect to $\mu$. Diverse proofs and generalizations of this contraction result are given in~\cite{otto05,wang-book,bgl2}. 

\medskip

A crucial challenging problem now consists in understanding the role of the dimension in the contraction property in Wasserstein distance. Indeed curvature and dimension are jointly considered in the synthetic definition by Lott-Sturm-Villani \cite{sturm,LV}, contraction properties, gradient commutation type properties or the Bakry-\'Emery curvature-dimension condition. It is for instance well known that, given $R\in\R$ and $n\geq1$, the $CD(R,n)$ curvature-dimension condition proposed by  D. Bakry and M. \'Emery in~\cite{bakryemery}, see section~\ref{sec-cd}, is satisfied for the Laplace-Beltrami operator on an $n$-Riemannian manifold if the Ricci curvature of the manifold is uniformly bounded from below by $R$.

This has been very recently performed in the following two remarkable results, by deriving an upper bound on the distance $W_2(H_t f\mu,H_sg\mu)$ with two {\it different} times $s,t>0$:

\medskip

$\bullet$ The first result is due to K. Kuwada in~\cite{kuwada} :  the Ricci curvature of the $n$-dimensional manifold $M$ is bounded from below by a constant $R\in\R$ if and only if 
\begin{equation}
\label{eq-contraction-kuwada}
W_2^2(H_t f\mu,H_sg\mu)\leq A(s,t,R) W_2^2( f\mu,g\mu)+B(s,t,n,R),
\end{equation}
for any $s,t>0$ and any probability densities $f,g$ with respect to $\mu$, for appropriate functions $A,B\geq0$. In the case $R=0$ the bound simplifies into
\begin{equation}
\label{eq-contraction-eks}
W_2^2(H_t f \mu,H_sg \mu)\leq W_2^2( f \mu,g \mu)+2 n(\sqrt{t}-\sqrt{s})^2,
\end{equation}
stated independently in~\cite{bgl2}.

\medskip

$\bullet$ The second result is due to M. Erbar, K. Kuwada and K.-T. Sturm~\cite{EKS13} : the Ricci curvature of the $n$-dimensional manifold $M$ is bounded from below by a constant $R\in\R$ if and only if
$$
s_{\frac Rn}\left(\frac 12 W_2(H_tf\mu,H_sg\mu)\right)^2\leq e^{-R(t+s)}\,s_{\frac Rn}\left(\frac 12 W_2(f\mu,g\mu)\right)^2+\frac nR(1-e^{-R(s+t)})\frac{(\sqrt{t}-\sqrt{s})^2}{2(t+s)}
$$
for any $s,t>0$ and any probability densities $f,g$ with respect to $\mu$. Here $s_r(x)=\sin(\sqrt{r}x) / \sqrt{r}$ if $r>0$, $s_r(x)=\sinh(\sqrt{-r}x) / \sqrt{-r}$ if $r<0$ and $s_0(x)=x$, hence recovering \eqref{eq-contraction-eks} for $R=0$.

\medskip

Observe that in these two results the dimension dependent additional term in the right-hand side is positive, and appears only when the two solutions are considered at different times $s$ and $t$. 

\medskip

A first aim of this paper is to take the dimension into account and to improve inequality~\eqref{ine-vonrenesse} for solutions considered at the {\it same} time. For instance in section~\ref{sec-first} we prove that 
\begin{equation}
\label{eq-contraction-un}
W_2^2(H_t fdx,H_tgdx)\leq W_2^2(fdx,gdx)-\frac{2}{n} \int_0^t\big(\ent{dx}{H_uf}-\ent{dx}{H_ug}\big)^2du
\end{equation}
for the heat semigroup on $\R^n$, any $t\geq0$ and any probability densities $f,g$ with respect to the Lebesgue measure $dx$; here $\ent{\mu}{g}=\int g\log gd\mu$ is the entropy. This inequality improves on~\eqref{ine-vonrenesse} since the Euclidean space $\R^n$ has null Ricci curvature and then satisfies~\eqref{ine-vonrenesse} with $R=0$. Let us observe that a dimensional contraction property in a Wasserstein distance with a modified cost was derived by Wang in~\cite{wang-cdrn}.

\medskip

A second aim it to obtain dimensional contraction inequalities for more general Markov semigroups. For that purpose we will work with a new distance called {\it Markov transportation distance}, based on the generator of the semigroup, and adapted to it and to the Bakry-\'Emery curvature-dimension condition formulation. 
It is defined by a modification of the following dynamical formulation of the Wasserstein distance proposed by J.-D. Benamou and Y. Brenier in~\cite{benamoubrenier} : for any probability densities $f$ and $g$ with respect to the Lebesgue measure in $\R^n$,  
$$
W_2(fdx,gdx)=\inf \Big( \int_0^1\int \frac{| w_s|^2}{\rho_s}\, dx ds \Big)^{1/2}
$$
where the infimum runs over all paths $(\rho_s)_{s\in[0,1]}$ and vector fields $(w_s)_{s\in[0,1]}$ such that $\partial_s\rho_s+\nabla\cdot w_s=0$, $\rho_{0}=f$ and $\rho_{1}=g$; here $\nabla\cdot$ stands for the divergence operator on $\rr^n$.  This dynamical approach is the starting point of the definition in~\cite{dns1,dns2} of generalized distances.  Instead, one can consider the quantity
$$
\inf \Big( \int_0^1\int \frac{|\nabla h_s|^2}{\rho_s} \, dx ds  \Big)^{1/2}
$$
where the infimum runs over all paths $(\rho_s, h_s)_{s\in[0,1]}$ such that $\partial_s \rho_s +\nabla\cdot (\nabla h_s) =0$, $\rho_0=f$ and $\rho_1 =g$.
This quantity is more adapted to the context of  a general Markov semigroup, and gives us the way to define the Markov transportation distance. Given a Markov generator $L$ on a space $E$, with carr\'e du champ $\Gamma$ and invariant measure $\mu$, we let
$$
T_2(f\mu,g\mu)=\inf \Big( \int_0^1\int \frac{\Gamma( h_s)}{\rho_s}d\mu ds \Big)^{1/2}
$$
for two probability densities $f$ and $g$ with respect to $\mu$, under the constraints $\partial_s  \rho_s+Lh_s=0, \rho_0 = f, \rho_1 = g$. In this abstract formulation, discrete and non-local operators can be studied in a similar way. A fundamental instance is that of $L = \Delta_g - \nabla V\cdot  \nabla$ on a Riemannian manifold $(E,g)$, with carr\'e du champ $\Gamma (f) = \vert \nabla f \vert^2$ and invariant measure $\mu$ with density $e^{-V}$. This will be the main example in this article, and in this case $W_2(f\mu,g\mu)\leq T_2(f\mu,g\mu)$ since the infima defining the distances run over a smaller set for $T_2$ than for $W_2$.

\medskip

The paper is organized as follows. In section~\ref{sec-first} we show in a simple way how to reach the dimension dependent contraction property~\eqref{eq-contraction-un} in Wasserstein distance for the specific heat semigroup on $\rr^n$. It will give a flavor of the results proved and the methods used below in the Markov transportation distance $T_2$ in our context of Markov semigroups $(P_t)_{t \geq 0}$ on a connected and compact Riemannian manifold.  

This distance is properly defined in section~\ref{sec-defusu}, together with fundamental properties and examples. In particular we derive an Otto-Villani theorem for $T_2$ : a logarithmic Sobolev inequality implies a transportation Talagrand inequality. 

Section~\ref{sec-cd} is devoted to our main application : the contraction property under the curvature-dimension condition $CD(R,n)$ on our semigroup. Under this condition we prove that 
$$
T_2^2(P_T f\mu,P_T g\mu)\leq e^{-2R T}T_2^2(f\mu,g\mu)-\frac{2}{n}\int_0^T e^{-2R(T-t)}\PAR{\ent{\mu}{P_t g}-\ent{\mu}{P_t f}}^2dt,
$$
for any $T>0$ and any probability densities $f,g$ with respect to the invariant measure $\mu$.

In section~\ref{sec-EVI} we briefly consider the so-called  evolution variational inequalities (EVI in short). These inequalities say that if the Ricci curvature of a manifold is bounded from below by a constant $R\in\R$, then\begin{equation}
\label{ine-evi}
W_2^2(f\mu, H_t g\mu)-W_2^2(f\mu,g\mu)\leq -\frac{e^{-2R t}-1+2R t}{2R t}W_2^2(f\mu,g\mu)+2t(\ent{\mu}{f}-\ent{\mu}{H_t g}).
\end{equation}
for the heat semigroup $(H_t)_{t \geq 0}$, any $t\geq0$ and any probability densities $f,g$ with respect to $\mu$. This inequality characterizes $(H_t)_{t \geq 0}$ as the gradient flow of the entropy with respect to the Wasserstein distance. This interpretation has been made by  R. Jordan, D. Kinderlehrer and F. Otto in~\cite{jko}, and has led to numerous developments, see in particular the seminal paper~\cite{ov00} and the huge contribution of~\cite{ambrosio-gigli-savare}. In section~\ref{sec-EVI} we explain how to obtain a dimensional EVI for the Wasserstein distance and the Euclidean heat semigroup, and then for  the Markov transportation distance and our Markov semigroup under a curvature-dimension condition. 

In Section~\ref{sec-gene} we briefly  investigate natural generalizations of the Markov transportation distance. 

\medskip

Many questions are left aside in this work, such as the general existence of geodesics, dual formulations and further equivalence between the obtained contraction and curvature conditions.  The purpose of this work is rather to show the interest of the $T_2$ distance, and these questions will be further investigated elsewhere. We haven chosen to present the Markov transportation distance in the classical setting of a compact Riemannian manifold to  properly prove the dimensional contraction inequality. We are convinced that the new distance can be defined in a general setting.  

 Since this work was completed, the second author \cite{gentil} has extended the bound~\eqref{eq-contraction-un} to the heat semigroup on an $n$-dimensional compact Riemannian manifold. Moreover, related results are studied in the work~\cite{AMS14} in preparation.

\section{The heat equation on $\R^n$}
\label{sec-first}

This section is devoted to the simple derivation of a dimension dependent contraction property for the heat semigroup $(H_t)_{t \geq 0}$ on $\rr^n$. It is defined by
$$
H_{t}f (x) = \int_{\rr^n} f(y) \frac{e^{-\frac{\vert x-y \vert^2}{4t}}}{\PAR{4\pi t}^{n/2}}\, dy
$$ 
and is obtained as the solution of the heat equation 
$\partial_t u=\Delta u;$
here $\Delta$ is the usual Laplace operator in $\R^n$.

For this semigroup, the bound \eqref{ine-vonrenesse} is classical with $R=0$ and $\mu$ the Lebesgue measure $dx$ on $\rr^n$, and is optimal in the sense that equality holds for all $t$ if $g$ is obtained from $f$ by a translation in $\rr^n$. Let us see how to simply turn this classical bound into a more precise dimension dependent bound.


Following~\cite{dns2}, let $(R_t)_{t\geq0}$ be the heat semigroup acting on $\rr^n$-valued maps, coordinate by coordinate. It satisfies 
\begin{equation}
\label{eq-dol}
H_t(\nabla\cdot w)=\nabla\cdot (R_t w)
\end{equation}
for all $\R^n$-valued functions $w$.  This semigroup acting on vectors will be the main tool in our derivation. We omit regularity issues which are carefully considered in~\cite{dns2}. 


As recalled in the introduction, the Benamou-Brenier Theorem ensures that 
\begin{equation}\label{eq-BB}
W_2^2(fdx,gdx)=\inf{\int_0^1\int \frac{|w_s|^2}{\rho_s}dsdx}
\end{equation}
for any probability measures $fdx$ and $gdx$ in $\R^n$; here the infimum runs over all couples $(\rho_s,w_s)_{s\in[0,1]}$ such that 
\begin{equation}
\label{eq-contrainte}
\partial_s \rho_s+\nabla\cdot w_s=0
\end{equation}
where, for all $s\in[0,1]$, $\rho_s$ is a probability density with respect to Lebesgue measure, $\rho_0=f$ and $\rho_1=g$.

Let now  $(\rho_s,w_s)_{s\in[0,1]}$ interpolate the densities $f$ and $g$ with the constraint~\eqref{eq-contrainte}. Then  $(H_t(\rho_s))_{s\in[0,1]}$ interpolates the densities $H_t f$ and $H_t g$ and, by~\eqref{eq-dol}, the couple $(H_t(\rho_s),R_t(w_s))_{s\in[0,1]}$ satisfies~\eqref{eq-contrainte}.  Then, by \eqref{eq-BB}, 
\begin{equation}
\label{eq-presqueder}
W_2^2(H_T fdx,H_Tgdx)\leq \int_0^1\int \frac{|R_T(w_s)|^2}{H_T(\rho_s)}dsdx
\end{equation}
for any $T\geq0$. Moreover:

\begin{elem}
\label{lem1-lambda}
Let $F:\dR^n \to\dR^n$  and $g$ a positive probability density with respect to the Lebesgue measure, with $F$ and $g$ smooth.  Then, for all $T \geq 0$
$$
\int\frac{|R_T F|^2}{H_T g}dx \leq \int\frac{|F|^2}{g}dx - \frac{2}{n}\int_0^T\PAR{\int \frac{R_t F\cdot \nabla H_tg}{H_t g}dx}^2dt.
$$
\end{elem}

\begin{eproof}
We let 
$$
\Lambda(t)=\int\frac{|R_t F|^2}{H_t g}dx
$$
for $t  \geq 0$ and prove that
$$
\Lambda'(t)\leq -\frac{2}{n}\PAR{\int \frac{R_t F\cdot \nabla H_t g}{H_t g}dx}^2,
$$
which will prove the lemma by time integration. Indeed 
$$
\Lambda'(t)=\int  \Big(2\frac{R_t F\cdot \Delta R_tF}{H_tg}-\frac{\Delta H_tg|R_t F|^2}{(H_tg)^2}\Big)dx.
$$
For notational simplicity, we let $\bar F=R_t F$, $\bar g=H_t g$ and then $\bar G=\log \bar g$. Since 
$$
0=\int \Delta\Big(\frac{|\bar F|^2}{\bar g}\big)dx=\int  2\nabla(|\bar F|^2) \cdot \nabla \Big(\frac{1}{\bar g}\Big)+\frac{1}{\bar g} \, \Delta (|\bar F|^2)+|\bar F|^2 \, \Delta \Big(\frac{1}{\bar g}\Big) \, dx,
$$
we obtain 
\begin{eqnarray*}
\Lambda'(t)
&=&
-\int\frac{2}{\bar g}\Big(\frac12\Delta|\bar F|^2-\bar F\cdot \Delta \bar F+\nabla (|\bar F|^2)\nabla \bar G+|\bar F|^2|\nabla \bar G|^2\Big)dx \\
&=&
-\int\frac{2}{\bar g}\sum_{1\leq i,j\leq n}\Big(\partial_i \bar F_i+\bar F_i\partial_j \bar G)^2dx 
\leq  -\int\frac{2}{\bar g}\sum_{1\leq i\leq n}\Big(\partial_i \bar F_i+\bar F_i\partial_i \bar G\Big)^2dx\\
&=&
 -\frac{2}{n}\int\bar g \, \Big(\sum_{1\leq i\leq n} \frac{\partial_i \bar F_i}{\bar g}+\frac{\bar F_i\partial_i \bar G}{\bar g} \Big)^2 \, dx\\
&\leq&
 -\frac{2}{n}\Big(\sum_{1\leq i\leq n}\int \big({\partial_i \bar F_i}+\frac{\bar F_i\partial_i \bar g}{\bar g}\big)dx\Big)^2 = -\frac{2}{n} \Big( \int \frac{\bar F \cdot \nabla \bar g}{\bar g} dx\Big)^2
\end{eqnarray*}
by the Cauchy-Schwarz inequality, the Jensen inequality for the probability measure  $\bar g \, dx$ and the relation $\displaystyle \int \sum_i \partial_i \bar F_i \, dx=0$.\end{eproof}

Then, by Lemma~\ref{lem1-lambda} and the Cauchy-Schwarz inequality (with respect to the measure $ds$), 
$$
\int_0^1\int \frac{|R_T(w_s)|^2}{H_T(\rho_s)}dsdx\leq \int_0^1\int \frac{|w_s|^2}{\rho_s}dsdx-\frac{2}{n}\int_0^T\PAR{\int \int_0^1 \frac{R_t (w_s)\cdot \nabla H_t(\rho_s)}{H_t (\rho_s)}dxds}^2dt. 
$$
Moreover the couple $(H_t(\rho_s),R_t(w_s))$ satisfies~\eqref{eq-contrainte}, so
$$
\int \frac{ R_t (w_s)\cdot \nabla H_t(\rho_s)}{H_t (\rho_s)}dx=\partial_s \int H_t(\rho_s)\log H_t(\rho_s)dx,
$$
and then 
$$
\int \int_0^1 \frac{R_t (w_s)\cdot \nabla H_t(\rho_s)}{H_t (\rho_s)}dxds = \ent{dx}{H_tf}-\ent{dx}{H_tg}.
$$

Then inequality~\eqref{eq-presqueder}  leads to the following refined contraction inequality for the heat semigroup in $\R^n$:

\begin{eprop}
\label{prop-facile1}
Let $(H_t)_{t\geq0}$ be the heat semigroup on $\dR^n$. Then for any  probability densities  $f$ and $g$ in $\dR^n$ such that $W_2(fdx,gdx)<\infty$,  for any $T>0$, 
\begin{equation}
\label{eq-contraction-heat}
W_2^2(H_T fdx,H_Tgdx)\leq W_2^2(fdx,gdx)-\frac{2}{n} \int_0^T\big(\ent{dx}{H_tf}-\ent{dx}{H_tg}\big)^2dt.
\end{equation}
\end{eprop}

\begin{erem} By comparison with~\eqref{eq-contraction-eks}, the dimension brings a negative correction term in the contraction property. The bound ~\eqref{eq-contraction-heat} is again an equality if $g$ is obtained from $f$ by translation. Finally, a Taylor expansion of~\eqref{eq-contraction-heat}, for $T$ close to $0$ and $g$ close to $f$, for any given $f$, implies back the curvature dimension $CD(0,n)$ for the Laplace operator (see section~\ref{sec-cd} below for the precise definition of the curvature-dimension condition). \end{erem}

\begin{erem}\label{contts}
Note that this result not only gives a correction for equal times, but also for different times $s,t$ in the spirit of \cite{bgl2}, \cite{EKS13} or \cite{kuwada}. Let indeed $s \leq t$: then applying the contraction estimate \eqref{eq-contraction-heat} to $P_{t-s}f$ and $g$ and then using \eqref{eq-contraction-eks} lead to
$$
W_2^2(H_t fdx,H_sgdx)\leq W_2^2(fdx,gdx)+n(t-s)-\frac{2}{n} \int_0^s\big(\ent{dx}{H_{t-s+u}f}-\ent{dx}{H_ug}\big)^2du.
$$
\end{erem}

Our main goal is then the extension of this contraction result to general Markov semigroups satisfying a $CD(R,n)$ condition, which will be given in Theorem~\ref{thm-contraction-rn} : there the Markov transportation distance will prove to be an adapted and efficient tool.

\section{The Markov transportation distance}
\label{sec-defusu}
\subsection{Definition}
\label{sec-defdef}

Let $(E, g)$ be a $C^{\infty}$ compact connected Riemannian manifold with volume measure $dx$, and $V$ be a smooth function on $E$ with $\int e^{-V} dx =1$. Let also $\mu$ be the Borel probability measure on $E$ with density $e^{-V}$. Let $(P_t)_{t\geq0}$ be the Markov semigroup on $E$ with infinitesimal generator 
$$
L = \Delta_g - \nabla V \cdot \nabla
$$ where $\Delta_g$ is the Laplace-Beltrami operator on $E$. The generator is defined on a dense subspace $\mathcal D(L)$ of $L^2(\mu)$, which includes the algebra $\mathcal A$ of (bounded) smooth functions (in our context, smooth means $C^\infty$). This algebra is itself
stable by $L$ and $P_t$. 

The semigroup is {\it reversible} with respect to $\mu$, in the sense that
\begin{equation}\label{eq-rev}
\int f \, P_t g \, d\mu = \int g \, P_t f \, d\mu \qquad {\textrm{or equivalently}} \qquad \int f \, L g \, d\mu = \int g \, L f \, d\mu
\end{equation}
for all $f, g \in \mathcal A$ and $t \geq 0$,

Such a Markov semigroup admits a Markov probability kernel, that is for any function $f\in L^2(\mu)$, $t\geq0$ and $x\in E$,   
$$
P_{t} f(x) = \int_{\rr^n} f(y) \, p_t(x, dy).
$$
The carr\'e du champ operator, defined on functions $f,g\in\mathcal A$ by the general expression
$$
\Gamma(f,g) = \frac{1}{2} \Big( L(fg) - f \, Lg - g \, Lf \Big) \in \mathcal A,
$$
satisfies $ \Gamma (f,f) = |\nabla f|^2$ where $|\nabla f|$ stands for the length of the vector $\nabla f$; for simplicity we shall let $\Gamma(f) = \Gamma(f,f)$. The Dirichlet form $\mathcal E_\mu(f)=\int \Gamma(f) d\mu$ is defined on its domain $\mathcal D(\mathcal E_\mu) \subset L^2(\mu)$, which also includes $\mathcal A$.

The generator satisfies the following so-called {\it diffusion} property :    
\begin{equation}
\label{eq-diffusion}
L\Phi(g)=\Phi'(g)Lg+\Phi''(g)\Gamma(g)
\end{equation}
for any smooth function $\Phi$ and any function $g\in\mathcal A$. In particular $\Gamma(\Phi'(g))=\Phi''^2(g)\Gamma(g)$.

We refer to~\cite{bgl1} for further details on these notions, and for their definition in a more general setting (there, the triple $(E, \Gamma, \mu)$ is called a compact Markov triple).

\medskip

Before giving the definition of the Markov transportation distance, we need to define the paths between probability densities.
We let $\mathcal F$ be the set of positive probability densities (with respect to $\mu$) in $\mathcal A$.

\begin{edefi}\label{defiT2}
 For a couple $(\rho_s,h_s)_{s\in[0,1]}$ of smooth functions on $[0,1] \times E$ with $\rho_s$ in $\mathcal F$ we define 
$$
\phi(\rho_s,h_s)=\int\frac{\Gamma(h_s)}{\rho_s}d\mu\in[0,\infty)
$$
for $s\in[0,1]$, and the action
$$
\Phi(\rho,h)=\int_0^1\phi(\rho_s,h_s)ds=\int_0^1\int\frac{\Gamma(h_s)}{\rho_s}d\mu ds\in[0,\infty].
$$

 For $f,g$ in $\mathcal F$ we call admissible path between $f$ and $g$ such a couple $(\rho_s,h_s)_{s\in[0,1]}$ for which moreover \begin{equation}
\label{eq-weak-sol}
\partial_s\rho_s+Lh_s=0, 
\end{equation}
$\rho_0=f$ and $\rho_1=g$. We let $\mathcal A(f,g)$ be the set of admissible paths between $f$ and $g$. 
\end{edefi}

\begin{edefi}\label{defiT2suite}
The Markov transportation distance is defined for $f,g\in\mathcal F$ by 
$$
T_2(f\mu,g\mu)=\inf \left( \int_0^1 \int \frac{\Gamma(h_s)}{\rho_s}d\mu ds \right)^{1/2},
$$ 
where the infimum runs over all admissible paths $(\rho_s,h_s)_{s \in [0,1]} \in\mathcal A(f,g)$. 
\end{edefi}

For given $f,g\in\mathcal F$, there exists a smooth function $h\in\mathcal A$  such that  $Lh=f-g$. In particular $\mathcal A(f,g)$ is nonempty since  $\rho_s=sf+(1-s)g$, associated with the function $h_s=h$ (independent of $s$), is an admissible  path. 

Moreover, $T_2 (f\mu, g\mu)$ is well defined and finite since  it is bounded from above by
$$
\int_0^1\!\!\phi(sf+(1-s)g,h) \, ds =  \int_0^1\!\!\int\! \frac{\Gamma(h)}{sf+(1-s)g} \; ds\,d\mu = \int {\Gamma(h)} \, \frac{\log(f)-\log(g)}{f-g} \, d\mu
$$
which is finite since $h\in\mathcal A$. For instance, if $f,g>\eta$ for some $\eta>0$, then
$$
T_2^2(f\mu,g\mu)\leq \frac1\eta\int {\Gamma(h)}d\mu<\infty. 
$$

\subsection{Remarks and examples}
\label{sec-exemple}

The Markov transportation distance heavily depends on both the reference measure $\mu$ and the generator $L$. 

\medskip

As it has been presented in the introduction, the Markov transportation distance is a generalization of the Benamou-Brenier dynamical formulation of the Wasserstein distance $W_2$ (see \cite{benamoubrenier}). In our setting of  a compact Riemannian manifold $E$ equipped with the probability measure $d\mu (x) = e^{-V(x)} dx$, it was indeed proven by F. Otto and M. Westdickenberg \cite[Prop. 4.3 and (4.15)]{otto05} that, given two positive densities (with respect to $\mu$) $f$ and $g$ in $\mathcal F$,
$$
W_2^2 (f \mu, g \mu) = \inf \int_0^1 \int \frac{\vert w_s \vert^2}{\rho_s} d\mu ds
$$
where the infimum runs over all smooth vector fields $(\rho_s, w_s)_{s \in [0,1]}$ on $[0,1] \times E$ with $\rho_s>0$ in $\mathcal F$ satisfying
\begin{equation}\label{eqcont}
\partial_s \rho_s + \nabla \cdot w_s - \nabla V \cdot w_s =0, \qquad \rho_0 = f, \rho_1 = g.
\end{equation}
By comparison, in this setting, 
$$
T_2^2(f\mu, g\mu) = \inf \int_0^1 \int \frac{\vert \nabla h_s\vert^2}{\rho_s} d\mu ds
$$
where the infimum runs over all smooth functions $(\rho_s, h_s)_{s \in [0,1]}$ satisfying
\begin{equation}\label{eqgenh}
\partial_s \rho_s + \Delta_g h_s - \nabla V \cdot \nabla h_s =0, \qquad \rho_0 = f, \rho_1 = g.
\end{equation}

But $(\rho_s, w_s = \nabla h_s)$ satisfies \eqref{eqcont} for any $(\rho_s, h_s)_{s \in [0,1]} \in \mathcal A (f,g)$, so
$T_2^2(f\mu, g\mu) \geq W_2^2 (f\mu, g\mu)$.

\bigskip

 The Markov transportation distance as defined in Definitions \ref{defiT2} and \ref{defiT2suite} can also be considered in the more general setting of a Polish measure space, a generator $L$, its carr\'e du champ $\Gamma$ and a reference measure $\mu$, but below we prefer to stick to our framework to be able to  properly justify our computation. Here is however an example in the discrete case.

\smallskip

In the case of a countable state space $E$, a Markov semigroup $(P_t)_{t \geq 0}$ is described by an infinite matrix of positive kernels ${(p_t(x,y))}_{(x,y)\in E \times E}$, $ t \geq 0$, such that for all $t \geq 0$
and $x \in E$, and any positive function $f$ on $E$,
$$
P_t f (x) = \sum_{y \in E}  f(y) \, p_t(x,y).
$$
For any $x\in E$, $p_t(x,.)$ is a probability measure on $E$. The generator $L$ is given by an infinite
matrix ${(L(x,y))}_{(x,y) \in E \times E}$, where for any finitely supported function $f$ on $E$,
$$
L f (x)= \sum_{y \in E} L(x,y) f(y).
$$
For the matrix $L$ to be a generator, it is required that $L (x,y) \geq 0$ whenever $x \neq y$, and $ \sum_{y } L(x,y)=0$ for every $x\in E$. The carr\'e du champ operator is defined on finitely  supported functions $f$ by 
$$
\Gamma(f)(x) = \frac{1}{2}\sum_{y \in E} L(x,y) \big [f(x)-f(y) \big ]^2, \quad x \in E.
$$
The measure is $\mu $ reversible if
$$
\mu(x) L(x,y)= \mu(y)  L(y,x). 
$$

\smallskip

Let us illustrate the discrete setting with the two point space $\{a,b\}$.  The generator is  unique up to a multiplicative factor, and is given by $Lf(a)=\kappa(f(b)-f(a))$ and $Lf(b)=\kappa( f(a)-f(b))$ for a nonnegative constant $\kappa$; moreover the carr\'e du champ is constant, equal to
$$
\Gamma(f)=\frac{\kappa}{2}(f(b)-f(a))^2,
$$
and the reversible measure is $\mu=\frac{1}{2}(\delta_a+\delta_b)$. There one can simply and explicitly compute a geodesic curve for the $T_2$ distance between two generic measures $(1-r) \delta_a + r \delta_b$ and $(1-t) \delta_a + t \delta_b$ with $0 < r, t < 1$.

\medskip

Let indeed $(\rho_s, h_s)_{s \in [0,1]} $ be an admissible path between $2 (1-r) \1_a + 2r \1_b$ and $2 (1-t) \1_a + 2t \1_b$. Then there exists a map $\phi:[0,1] \to [0,1]$ such that $\phi( 0)=r$ and $\phi(1)=t$, and 
$$
\rho_s=2 \phi(s)\1_b+2 (1-\phi(s))\1_a.
$$
The map $h_s$ has to satisfy $2 \phi'(s)(\1_b-\1_a)= - Lh_s$ for $s\in[0,1]$, that is $\phi'(s)=(h_s(b)-h_s(a))\kappa/2$. It remains to minimize 
$$
\int_0^1\int \frac{\Gamma(h_s)}{\rho_s}d\mu ds= \frac{\kappa}{2}  \int_0^1\int \frac{(h_s(b)-h_s(a))^2}{\rho_s}d\mu ds= \frac{1}{\kappa} \int_0^1\phi'(s)^2\Big(\frac{1}{\rho_s(a)}+\frac{1}{\rho_s(b)}\Big)ds
$$
Since $\rho_s(a)=2(1-\phi(s))$ and $\rho_s(b)=2 \phi(s)$ we need to minimize
$$
\frac{1}{2 \kappa}\int_0^1\phi'(s)^2\Big(\frac{1}{\phi(s)}+\frac{1}{1-\phi(s)}\Big)ds
$$
over all functions $\phi$ such that $\phi(0)=r$ and $\phi(1)=t$.  The Euler-Lagrange equation is
$$
2\phi''\Big(\frac{1}{\phi}+\frac{1}{1-\phi}\Big)=\phi'^2\Big(\frac{1}{\phi^2}-\frac{1}{(1-\phi)^2}\Big).
$$
It implies that ${\phi'^2}=a\phi(1-\phi)$ for some $a>0$. We let $r = \sin^2 \theta$ and $t = \sin^2 \omega$ with $0 < \theta , \omega < \pi/2$. This solves into $\phi(s)=\sin^2(\frac{\sqrt{a}}{2}s + \varepsilon \theta)$ for $s\in[0,1]$, and some $\varepsilon = \pm 1$. For such a $\varphi$,
$$
\phi'(s)^2\Big(\frac{1}{\phi(s)}+\frac{1}{1-\phi(s)}\Big)= \frac{\phi'^2}{\phi(1-\phi)}=a, \quad s\in[0,1].
$$
But the smallest $a$ for which $\varphi(1) = t$ is $a= 4 (\omega - \theta)^2$ (and $\varepsilon = 1$) : hence $\varphi$ is given by $\varphi (s) = \sin^2(s \omega + (1-s) \theta)$. This implies that 
$$
T_2^2\big( (2 (1-r) \1_a + 2r \1_b)\mu,(2 (1-t) \1_a + 2t \1_b)\mu \big)= \frac{2 (\omega - \theta)^2}{ \kappa} 
$$ 
for the $T_2$ distance defined as in Definitions \ref{defiT2} and \ref{defiT2suite}. Moreover 
$$
\int \frac{\Gamma(h_s)}{\rho_s}d\mu= \frac{2 (\omega - \theta)^2}{ \kappa}
$$
for all $s\in[0,1]$, so for such a $\varphi$ the path $(\rho_s, h_s)$ is a geodesic between $2 (1-r) \1_a + 2r \1_b$ and $2 (1-t) \1_a + 2t \1_b$.

\bigskip

 In this example, in a discrete setting, one has obtained the existence of a  geodesic path between two probability measures associated to the Markov transportation distance. Let us recall that there is no such geodesic for the Wasserstein distance.

\subsection{General properties of $T_2$}

We first exhibit $\ep$-geodesics for the $T_2$ distance.  Actually, we shall see below how properties on the distance and curvature-dimension bounds can be obtained without geodesics.

\begin{eprop}[$\ep$-geodesics]
\label{prop-geo}
Let $f,g\in\mathcal F$ and let $\ep>0$.  Then there exists an $\ep$-geodesic map, that is an admissible path  $(\rho_s,h_s)_{s \in [0,1]} \in \mathcal A (f,g)$ such that for all $s\in[0,1]$, 
$$
\phi(\rho_s,h_s)=\int \frac{\Gamma(h_s)}{\rho_s}d\mu\leq T_2^2(f\mu,g\mu)+\ep.
$$ 
\end{eprop}

\begin{eproof}
Let $\ep>0$ and an admissible path  $(\rho_s,h_s)\in\mathcal A(f,g)$ such that 
$$
\Phi(\rho, h) =  \int_0^1\int \frac{\Gamma(h_s)}{\rho_s}d\mu ds \leq T_2^2(f\mu,g\mu)+\ep.
$$ 
It is easy to see that there exists $a>0$  such that 
$$
{\int_0^1\sqrt{\phi(\rho_u,h_u)+a} \,du}={\sqrt{\disp\Phi(\rho,h)+\ep}}.
$$
Then let $\beta: [0,1] \to [0,1]$ be defined by  
$$
s=\frac{\disp\int_0^{\beta(s)}\sqrt{\phi(\rho_u,h_u)+a} \,du}{\sqrt{\disp\Phi(\rho,h)+\ep}}
$$
for $s\in[0,1]$. The function $\beta$ is increasing and differentiable in $[0,1]$  and satisfies $\beta(0)=0$ and $\beta(1)=1$, so
 $(\rho_{\beta(s)},\beta'(s)h_{\beta(s)})\in\mathcal A(f,g)$. Moreover
 \begin{eqnarray*}
 \phi(\rho_{\beta(s)},\beta'(s)h_{\beta(s)})=\beta'(s)^2\phi(\rho_{\beta(s)},h_{\beta(s)})
 &=&
 \PAR{\Phi(\rho,h)+\ep}\frac{\phi(\rho_{\beta(s)},h_{\beta(s)})}{\phi(\rho_{\beta(s)},h_{\beta(s)})+a}\\
& \leq&
 \Phi(\rho,h)+\ep
 \; \leq \; T_2^2(f\mu,g\mu)+2\ep
 \end{eqnarray*}
for any $s\in[0,1]$. This means that the couple $(\rho_\beta,\beta'h_\beta)$ is a $2\ep$-geodesic. 
\end{eproof}

\begin{eprop}
\label{prop-demimetrique}
 The space $(\mathcal  F,T_2)$ is a metric space.  
\end{eprop}



\begin{eproof}
For any $f\in\mathcal F$, then $T_2(f\mu,f\mu)=0$ by choosing $\rho_s$ constant equal to $f$. Conversely, if $f$ and $g$ in $\mathcal F$ are such that $T_2(f\mu,g\mu)=0$, then $W_2(f\mu,g\mu)=0$ since it is smaller than $T_2(f\mu,g\mu)$, as seen in section \ref{sec-exemple}; hence $f = g$. Moreover $T_2$ is a symmetric function with respect to the two densities. 

Let now $f$, $g$ and $h$ in $\mathcal F$.  Let $(\rho^1_s, h^1_s)$ (resp. $(\rho^2_s, h^2_s)$) be an $\ep$-geodesic map between $f$ and $g$ (resp. $g$ and $h$). Let $\alpha\in(0,1)$ and define $(\rho_s, h_s)_{s \in [0,1]}$ by 
$$
\rho_s=
\left\{
\begin{array}{lr}
\disp\rho_{s/\alpha }^1,& \hbox{if }s\in[0,\alpha],\\
\disp\rho_{(s-\alpha)/(1-\alpha )}^2, &\hbox{if }s\in[\alpha,1].\\
\end{array}
\right.
\qquad 
h_s=
\left\{
\begin{array}{lr}
\disp \frac{1}{\alpha}h_{s/\alpha }^1,& \hbox{if }s\in[0,\alpha),\\
\disp \frac{1}{1-\alpha} h_{(s-\alpha)/(1-\alpha) }^2, &\hbox{if }s\in[\alpha,1].\\
\end{array}
\right.
$$
Then the couple $(\rho_s, h_s)$ is an admissible path between $f$ and $h$. Moreover
 $$
T_2^2(f \mu, h \mu) \leq \Phi(\rho,h)=\int_0^\alpha \! \PAR{\frac{1}{\alpha}}^2\phi(\rho_{s/\alpha }^1,h_{s/\alpha }^1)ds+\int_\alpha^1 \! \PAR{\frac{1}{1-\alpha}}^2\phi(\rho_{(s-\alpha)/(1-\alpha) }^2,h_{(s-\alpha)/(1-\alpha) }^2)ds.
$$
Now $\rho_s^1$ and $\rho_s^2$ are $\ep$-geodesics, so
$$
\phi(\rho_{s/\alpha }^1,h_{s/\alpha }^1)\leq T_2^2(f\mu,g\mu)+\ep
$$
 and 
 $$
 \phi(\rho_{(s-\alpha)/(1-\alpha) }^2,h_{(s-\alpha)/(1-\alpha) }^2)\leq T_2^2(g\mu,h\mu)+\ep.
 $$
Hence
$$
T_2^2(f \mu, h \mu) \leq {\frac{1}{\alpha}}(T_2^2(f\mu,g\mu)+\ep)+ {\frac{1}{1-\alpha}}(T_2^2(g\mu,h\mu)+\ep).
$$
Now, choose 
$$
\alpha=\frac{T_2(f\mu,g\mu)}{T_2(f\mu,g\mu)+T_2(g\mu,h\mu)}
$$ 
and let $\ep$ go to $0$ to obtain the triangular inequality
$$
T_2(f \mu, h \mu) \leq T_2(f\mu,g\mu)+T_2(g\mu,h\mu). 
$$
\end{eproof}

\begin{eprop}[Tensorization]
\label{prop-tensorisation} Let $(P_t^i)_{t\geq0}$, $i\in\{1,\cdots, N\}$ be $N$ Markov semigroups on compact connected Riemannian manifolds $E_i$ with probability measure  $\mu_i$,  with generators $L_i$ and carr\'es du champ $\Gamma_i$ as in Section~\ref{sec-defdef}.   Then one can define a product semigroup $P_t=\otimes_{i=1}^N P_t^i$ on the product space $E=\times_{i=1}^N E_i$ with $\mu=\otimes_{i=1}^N\mu_i$ with generator $L=\oplus_{i=1}^N L_i$ and carr\'e du champ $\Gamma=\oplus_{i=1}^N\Gamma_i$. 

Then, for any densities $f(x)=\prod_{i=1}^N f_i(x_i)$ and $g(x)=\prod_{i=1}^N g_i(x_i)$ ($x=(x_1,\cdots,x_N)$) in $\mathcal F$,
\begin{equation}
\label{eq-tenso1}
T_2^2(f\mu,g\mu)\geq \sum_{i=1}^NT_{2,i}^2(f_i\mu_i,g_i\mu_i).
\end{equation}
\end{eprop}

\begin{eproof}
For simplicity we  prove the result for $N=2$. Let $(\rho_s,h_s)$ be an admissible path between the densities $f_1(x)f_2(y)$ and $g_1(x)g_2(y)$. Let $\rho_s^1(x)=\int \rho_s(x,y)d\mu_2(y)$ and $h_s^1(x)=\int h_s(x,y)d\mu_2(y)$, and let $\rho_s^2$ and $h_s^2$ similarly defined. Then  
\begin{equation}\label{sommetenso}
\phi(\rho_s,h_s)=\int\frac{\Gamma(h_s)}{\rho_s}d\mu_1d\mu_2=\int\frac{\Gamma_1(h_s)}{\rho_s}d\mu_1d\mu_2+\int\frac{\Gamma_2(h_s)}{\rho_s}d\mu_1d\mu_2. 
\end{equation}

Let us first prove that 
\begin{equation}
\label{eq-jensen}
\int\frac{\Gamma_1(h_s)}{\rho_s}d\mu_2\geq \frac{\Gamma_1(h_s^1)}{\rho_s^1},
\end{equation}
and similarly for the second coordinate. 
Since 
\begin{equation}\label{eq-gamma}
\Gamma_1(f)(x)=\lim_{t\rightarrow0}\frac{1}{2t}\int\int(f(y_1)-f(y_2))^2p_t^1(x,dy_1)p_t^1(x,dy_2)
\end{equation}
for every function $f$, and for the Markov kernel $p_t^1$ of the semigroup $(P_t^1)_{t \geq 0}$ (see for instance~\cite{bgl1}), then for all $x$ 
\begin{multline*}
\int\frac{\Gamma_1(h_s)(x,y)}{\rho_s(x,y)}d\mu_2(y)=\\
\rho_s^1(x)\lim_{t\rightarrow0}\frac{1}{2t}\int \int\int\left(\frac{h_s(z_1,y)-h_s(z_2,y)}{\rho_s(x,y)}\right)^2 \frac{\rho_s(x,y)}{\rho_s^1(x)}\, d\mu_2(y) \, p_t^1(x,dz_1) \, p_t^1(x,dz_2)\\
\geq \frac{1}{\rho_s^1(x)} \lim_{t\rightarrow0}\frac{1}{2t}\int \int (\int (h_s(z_1,y)-h_s(z_2,y))d\mu_2(y))^2 \, p_t^1(x,dz_1) \, p_t^1(x,dz_2)
=\frac{\Gamma_1(h_s^1)(x)}{\rho_s^1(x)}
\end{multline*}
by the Cauchy-Schwarz inequality for the probability measure $\frac{\rho_s(x,y)}{\rho_s^1(x)}d\mu_2(y)$.  

By \eqref{sommetenso} and \eqref{eq-jensen} written for both variables we obtain 
$$
\phi(\rho_s,h_s)\geq \phi_1(\rho^1_s,h^1_s)+\phi_2(\rho^2_s,h^2_s). 
$$
After integration over $s\in[0,1]$, we get, for any admissible path $(\rho_s,h_s)\in\mathcal A(f_1f_2,g_1g_2)$~: 
\begin{equation}
\label{eq-tenso2}
\Phi(\rho,h)=\int_0^1\phi(\rho_s,h_s)ds\geq \Phi_1(\rho^1,h^1)+\Phi_2(\rho^2,h^2) \geq T_2^2(f_1\mu_1,g_1\mu_1)+T_2^2(f_2\mu_2,g_2\mu_2)
\end{equation}
since $(\rho_s^i,h_s^i)$  is an admissible path between  $f_i$ and $g_i$, for $i=1,2$. The result follows by optimizing over $(\rho_s, h_s)$. \end{eproof}

\subsection{First application :  the Talagrand inequality}

 As explained in the introduction, we will recover classical bounds as the contraction properties in the Markov transportation distance. We first make some observations on related functional inequalities.
\medskip

The so-called Otto-Villani Theorem says that a logarithmic Sobolev inequality with constant $C$ (see~\eqref{eqLS} below) implies the Talagrand inequality
$$
W_2^2(f \mu, \mu) \leq 4C \, \ent{\mu}{f}
$$
for all probability densities $f$. This inequality has first been derived by M. Talagrand in~\cite{talagrand96} and linked with the logarithmic Sobolev inequality in~\cite{ov00} (see also~\cite{bgl01}). The proofs in~\cite{ov00,gigliledoux13} are based on the general inequality
\begin{equation}
\label{eq-kuwadaw}
W_2^2(P_t f\mu, f\mu)\leq t (\ent{\mu}{f}-\ent{\mu}{P_t f})
\end{equation}
(see also~\cite{giglikuwada13}). We prove the same inequality for the larger $T_2$ distance : 

\begin{eprop}
\label{prop-kuwada}
For our diffusion semigroup $(P_t)_{t \geq 0}$ and any $f$ in $\mathcal F$, there holds
\begin{equation}
\label{eq-kuwada}
T_2^2(P_t f\mu, f\mu)\leq t \, (\ent{\mu}{f}-\ent{\mu}{P_t f})
\end{equation}
for every $t\geq0$, and in particular
\begin{equation}
\label{eq-derivative}
\limsup_{t\rightarrow 0^+}\frac{T_2(P_t f\mu, f\mu)}{t}\leq \sqrt{\int \frac{\Gamma(f)}{f}d\mu}.
\end{equation} 
\end{eprop}

\begin{eproof}
Let $f$ in $\mathcal F$ be given. Then $(P_{st}f, - tP_{st} f)_{s\in[0,1]}$ is an admissible couple between  $f$ and $P_t f$. By definition of the distance $T_2$, it implies that
$$
T_2^2(P_t f \mu, f\mu)\leq t\int_0^t\int \frac{\Gamma(P_r f)}{P_r f}d\mu dr
$$
 by change of time variable. Moreover  
$$
\frac{d}{dr} \ent{\mu}{P_r f} = \frac{d}{dr} \int P_r f\log P_r f \, d\mu =-\int \frac{\Gamma(P_rf)}{P_r f} \, d\mu
$$
by diffusion property of the semigroup. This leads to~\eqref{eq-kuwada} by integrating in $r$. 
\end{eproof}

\begin{ecor}\label{corologsob}
In our notation, assume that the probability measure $\mu$ satisfies a logarithmic Sobolev inequality with constant $C$, that is, 
\begin{equation}\label{eqLS}
\ent{\mu}{f}\leq C\int \frac{\Gamma(f)}{f}d\mu, 
\end{equation}
for any $f$ in $\mathcal F$. Then 
$$
T_2^2(f\mu, P_T f \mu)\leq 4C\ent{\mu}{f}
$$
for any $f$ in $\mathcal F$. In particular, if for instance $T_2$ is lower semicontinuous with respect to narrow convergence, then
$\mu$ satisfies a Talagrand type inequality for the distance $T_2$, namely
\begin{equation}\label{tala}
T_2^2(f\mu,\mu)\leq 4C\ent{\mu}{f}
\end{equation}
for any $f$ in $\mathcal F$.
\end{ecor}

\begin{eproof}
Let $f$ in $\mathcal F$ be given, and let $\phi(t)=\ent{\mu}{P_t f}$. Then~\eqref{eq-derivative} and semigroup properties imply  
$$
\frac{d^+}{dt}T_2(P_t f\mu, f\mu)\leq \sqrt{-\phi'(t)}.
$$
 Moreover the logarithmic Sobolev inequality for $\mu$ ensures that $\phi'(t)\leq -\phi(t) /C$, and thus 
$$
\sqrt{-\phi'(t)}\leq -\sqrt{{4C}}\PAR{\sqrt{\phi(t)}}'.
$$
Let now $T>0$. Then 
\begin{multline*}
T_2(P_T f\mu,f \mu)\leq \int_0^T\frac{d^+}{dt}T_2(P_t f\mu,f \mu)dt\leq -\sqrt{{4C}}\int_0^T \big( \sqrt{\phi(t)} \big)'dt\\
= \sqrt{{4C}}\PAR{\sqrt{\ent{\mu}{f}}-\sqrt{\ent{\mu}{P_Tf}}} \leq \sqrt{4C} \, \sqrt{\ent{\mu}{f}}.
\end{multline*}
Moreover $P_T f \mu$ narrowly converges to $f \mu$. Hence
$$
T_2(\mu, f\mu)\leq \liminf_{T\rightarrow\infty}T_2(P_Tf\mu,f\mu)\leq\sqrt{{4C}}\sqrt{\ent{\mu}{f}}
$$
if $T_2$ is lower semicontinuous with respect to narrow convergence.
\end{eproof}

Inequality \eqref{tala} will be useful in the following section in the derivation of refined convergence rates.

 
\section{Contraction property under the curvature-dimension condition $CD(R,n)$}
\label{sec-cd}
In this section we prove a dimension dependent contraction property in the Markov transportation distance. We will see that the $\Gamma_2$-calculus is a well adapted and efficient tool.

\subsection{Curvature condition, examples and useful commutation properties}

The $\Gamma_2$-operator, or iterated carr\'e du champ operator, is defined on functions $f\in\mathcal A$ by the general expression
 $$
 \Gamma_2(f) = \frac{1}{2} \Big( L \Gamma(f) - 2\Gamma (f, Lf) \Big).
 $$

\begin{edefi}
The  Markov semigroup $(P_t)_{t \geq 0}$ (associated to the generator $L$) is said to satisfy a curvature-dimension condition $CD(R,n)$ for $R\in\R$ and $n\geq 1$ if
$$
\Gamma_2(f)\geq \rho \, \Gamma(f)+\frac{1}{n}(Lf)^2,
$$
for all functions $f\in\mathcal A$. 
\end{edefi}
This criterion has been introduced in the seminal paper \cite{bakryemery} by D. Bakry and M. \'Emery. 
For instance the Laplace-Beltrami operator on the sphere $S^{n}\subset \R^{n+1}$ satisfies a $CD(n-1,n)$ condition.  
More generally, for a complete Riemannian manifold $M$ with dimension $d$, equipped with the Laplace-Beltrami operator $\Delta_g$ and the Riemannian measure $dx$, the curvature-dimension condition $CD(R,d)$ holds for $\Delta_g$ if the Ricci curvature of $M$ is uniformly bounded from below by $R$. An explicit condition on the potential $V$ which is equivalent to a $CD(R,n)$ condition for $L = \Delta_g - \nabla V \cdot \nabla$ is given in \cite[Appendix C.6]{bgl1} for instance. Observe that $n$ need not be the dimension on the manifold.

\medskip

 One of the main results concerning the curvature-dimension condition $CD(R,\infty)$ is a regularity property of the Markov semigroup. The $CD(R,\infty)$ conditions holds for a diffusive Markov semigroup if and only if for any function $f\in\mathcal A$
\begin{equation}
\label{eq-strong}
\Gamma(P_t(f))\leq e^{-2Rt}\PAR{P_t\sqrt{\Gamma(f)}}^2.
\end{equation}
This result, proved in~\cite{bakrysaintflour}, is the key point for many applications such as logarithmic Sobolev inequalities, Harnack parabolic inequalities, etc. (see~\cite{bgl1}). 
Gradient bounds, in a weaker form, also hold under the $CD(R,n)$ condition with finite $n$ (see~\cite{bakryledoux-liyau} and~\cite{wang-cdrn}). Here is   a new such bound which will be the key point for the Markov transportation distance. 

\begin{elem}
\label{lem4new}
Let $(P_t)_{t\geq0}$ be a diffusion Markov semigroup, and let $R\in\R$ and $n\geq 1$. The following assertions are equivalent :
\begin{enumerate}[(i)]
\item The Markov semigroup satisfies a $CD(R,n)$ condition. 
\item For all functions  $f,g\in\mathcal A $ with $g>0$ and all  $t \geq 0$
\begin{equation}
\label{eq-contracbis}
\frac{\Gamma(P_tf)}{P_t g}\leq e^{-2Rt}P_t\PAR{\frac{\Gamma(f)}{g}}-\frac{2}{n}\int_0^te^{-2Ru} \frac{\big[ LP_tf-P_u(\Gamma(P_{t-u} f,\log P_{t-u} g))\big]^2}{P_tg} du.
\end{equation}
\end{enumerate}

In particular, under the $CD(R,n)$ condition and for any $t\geq 0, f \in\mathcal A $ and $g \in \mathcal F$,
\begin{equation}
\label{eq-contracint}
\int \frac{\Gamma(P_tf)}{P_t g}d\mu\leq e^{-2Rt}\int \frac{\Gamma(f)}{g}d\mu-\frac{2}{n}\int_0^te^{-2R(t-u)}\PAR{\int \frac{\Gamma(P_u f,P_u g)}{P_ug}d\mu}^2du.
\end{equation}
 
\end{elem}
\begin{eproof}
Let us first prove that $(i)$ implies $(ii)$. We let $t>0$ and $f,g\in\mathcal A$ be fixed, with $g>0$. Then we define
$$
\Lambda(s)=P_s\PAR{\frac{\Gamma(P_{t-s}f)}{P_{t-s} g}}
$$
for $s\in[0,t]$, and then $F=P_{t-s}f$ and $G=P_{t-s}g$. Then 
$$
\Lambda'(s)=P_s\PAR{-2\frac{\Gamma(F,LF)}{G}+\Gamma(F)\frac{LG}{G^2}+L\PAR{\frac{\Gamma(F)}{G}}}.
$$
But
$$
L(hk)=2\Gamma(h,k)+hLk+kLh
$$
for any function $h,k\in\mathcal A$, so the diffusion property~\eqref{eq-diffusion} and the definition of $\Gamma_2$ lead to 
$$
\Lambda'(s)=2P_s\PAR{\frac{1}{G}\big[\Gamma_2(F)-\Gamma(\Gamma(F),\log G)+\Gamma(F)\Gamma(\log G)\big]} .
$$
Now Lemma~\ref{lem-technique} below, applied with $f=F$ and $g=-\log G$, ensures that  
$$
\Gamma_2(F)-\Gamma(\Gamma(F),\log G)+\Gamma(F)\Gamma(\log G)\geq R \, \Gamma(F)+\frac{1}{n}(LF-\Gamma(F,\log G))^2.
$$
Since $G\geq0$, this gives 
$$
\Lambda'(s) \geq 2R \Lambda(s)+\frac{2}{n}P_s\PAR{\frac{\big[LF-\Gamma(F,\log G)\big]^2}{G}} \geq 2R\Lambda(s)+\frac{2}{n}\frac{\big[LP_tf-P_s(\Gamma(P_{t-s} f,\log P_{t-s} g))\big]^2}{P_tg}
$$
by the Cauchy-Schwarz inequality for the Markov kernel of $P_s$ and semigroup properties. Inequality~\eqref{eq-contracbis} follows by integration over $s\in[0,t]$. 

\smallskip

Let us now assume $(ii)$ and let $g=1$. Then inequality~\eqref{eq-contracbis} writes 
$$
\Gamma(P_tf)\leq e^{-2Rt}P_t\Gamma(f)-\frac{2}{n}\int_0^te^{-2Ru}(LP_tf)^2du.
$$ 
Taking the time derivative at $t=0$ implies back the $CD(R,n)$ condition. 

\smallskip

Let us finally prove~\eqref{eq-contracint}:  integrating~\eqref{eq-contracbis} with respect to $\mu$ gives 
$$ 
\int \frac{\Gamma(P_tf)}{P_t g}d\mu\leq e^{-2Rt}\int \frac{\Gamma(f)}{g}d\mu-\frac{2}{n}\int_0^t e^{-2Ru} \Big[ \int \frac{\big[ LP_tf-P_u(\Gamma(P_{t-u} f,\log P_{t-u} g))\big]^2}{P_ug}d\mu \Big] du
$$
by invariance property of $\mu$. Then the Cauchy-Schwarz inequality for the measure $\mu$ implies~\eqref{eq-contracint} by recalling that $\int gd\mu=1$, invariance property of $\mu$ and change of time variable.
\end{eproof}

\begin{elem}
\label{lem-technique}
For a diffusion Markov semigroup, under the curvature-dimension condition $CD(R,n)$ (with $R\in\R$), for all functions $f,g\in\mathcal A$ 
\begin{equation}
\label{eq-technique}
\Gamma_2(f)+\Gamma(\Gamma(f), g)+\Gamma(f)\Gamma(g)\geq R\, \Gamma(f)+\frac{1}{n}\PAR{Lf+\Gamma(f,g)}^2.
\end{equation}
\end{elem}
\begin{eproof}
The proof is inspired from  Lemma~5.4.4, p.~83 of~\cite{logsob}. Let $f,g\in\mathcal A$ and $x_0\in E$. 
Let $\Phi$ be a smooth map on $\R^2$ such that 
$$
\partial_2\Phi=\partial^2_{11}\Phi=\partial^2_{22}\Phi=0,\quad\partial_1\Phi=1\quad {\rm and}\quad\partial^2_{12}\Phi=\frac{1}{2}
$$
at the point $(f(x_0),g(x_0))$.  Then the $CD(R,n)$ condition applied to the function $\Phi(f,g)$ at the point $x_0$ yields
$$
\Gamma_2(\Phi(f,g))\geq R \, \Gamma \Phi (f,g)+\frac{1}{n}(L\Phi(f,g))^2. 
$$
The usual change of variable rules for the $\Gamma$ and $\Gamma_2$ operators  (see for instance \cite[p. 83]{logsob}) imply 
$$
\Gamma_2(f)+\Gamma(\Gamma(f), g)+\frac{1}{2}\big[\Gamma(f,g)^2+\Gamma(f)\Gamma(g)\big]\geq R\Gamma(f)+\frac{1}{n}
\PAR{Lf+\Gamma(f,g)}^2. 
$$
The result follows since $\Gamma(f,g)^2\leq \Gamma(f)\Gamma(g)$.
\end{eproof}

\begin{erem}
Without the dimension, namely under the curvature-dimension condition $CD(R,\infty)$, inequality~\eqref{eq-contracbis}  is a direct consequence of inequality~\eqref{eq-strong}. Indeed~\eqref{eq-strong} implies  
$$
\frac{\Gamma(P_t f)}{P_t g}\leq e^{-2 R t}\frac{P_t \big(\sqrt{ \Gamma(f)} \, \big)^2}{P_t g}\leq e^{-2 R t}P_t\PAR{\frac{\Gamma(f)}{g}},
$$
by the Cauchy-Schwarz inequality for the Markov kernel of $P_t$.
\end{erem}

\subsection{Contraction property under $CD(R,n)$}

The following result is a extension of Proposition~\ref{prop-facile1} for the heat semigroup on $\rr^n$ to our diffusion Markov semigroup on a manifold under the curvature-dimension condition $CD(R,n)$.  For $R=0$ we precisely recover the bound obtained in Wasserstein distance for the heat semigroup on $\rr^n$.

\begin{ethm}
\label{thm-contraction-rn}
Let $(P_t)_{t \geq 0}$ be our diffusion Markov semigroup on $E$ satisfying a $CD(R,n)$ condition with $R\in\dR$ and $n\geq 1$. Then, for all $f,g\in\mathcal F$ and $T\geq0$, 
\beq
\label{eq-contraction-T2}
T_2^2(P_T f\mu,P_T g\mu)\leq e^{-2R T}T_2^2(f\mu,g\mu)-\frac{2}{n}\int_0^T e^{-2R(T-t)}\PAR{\ent{\mu}{P_t g}-\ent{\mu}{P_t f}}^2dt.
\eeq
\end{ethm}
\begin{eproof}
Let $(\rho_s,h_s)$ be an admissible  path between $f$ and $g$. Then $(P_t(\rho_s),P_t(h_s))_{s\in[0,1]}$ is also an  admissible path  between $P_tf$ and $P_t g$. 

Then inequality~\eqref{eq-contracint} of Lemma~\ref{lem4new}, applied at time $T$ to the functions $h_s$ and $\rho_s$, implies 
$$
\int\frac{\Gamma(P_T h_s)}{P_T\rho_s}d\mu \leq e^{-2R T}\int\frac{\Gamma( h_s)}{\rho_s}d\mu-\frac{2}{n}\int_0^Te^{-2R (T-t)}\PAR{\int \frac{\Gamma(P_t(\rho_s),P_t(h_s))}{P_t(\rho_s)}d\mu}^2dt.
$$
Integrating over $s\in[0,1]$ and the Cauchy-Schwarz inequality imply 
\begin{multline*}
\int_0^1\int\frac{\Gamma(P_T h_s)}{P_T\rho_s}d\mu\, ds\leq e^{-2R T}\int_0^1\int\frac{\Gamma( h_s)}{\rho_s}d\mu\, ds \\
-\frac{2}{n}\int_0^Te^{-2R (T-t)}\PAR{\int_0^1\int \frac{\Gamma(P_t(\rho_s),P_t(h_s))}{P_t(\rho_s)}d\mu ds}^2dt.
\end{multline*}
We finally obtain~\eqref{eq-contraction-T2} since, letting $\varphi(s) = \ent{\mu}{P_t \rho_s}$, 
$$
\int_0^1\int \frac{\Gamma(P_t(\rho_s),P_t(h_s))}{P_t(\rho_s)}d\mu \,ds=\int_0^1 \varphi'(s)\, ds = \varphi(1) - \varphi(0) = \ent{\mu}{P_t g}-\ent{\mu}{P_t f}
$$
by the reversibility property \eqref{eq-rev} of the semigroup.
\end{eproof}

\begin{erem}
As noted in the introduction, this result has a particular new flavor. Indeed the recent results \eqref{eq-contraction-kuwada}-\eqref{eq-contraction-eks} in \cite{bgl2}, \cite{EKS13} and \cite{kuwada} present a dimensional correction term for the contraction property, but for solutions at different times only. If the approaches for these inequalities are slightly different, it would be of interest to obtain a dimensional correction term for our contraction also in different times. A possible approach could be through Evolution variational inequalities, as studied in the next section, as the contraction result in \cite{EKS13} is deduced from these inequalities.
\end{erem}

Assuming that $\mu$ is a probability measure, and taking $g=1$, under the $CD(R,\infty)$ condition, the following bound
$$
T_2^2(P_T f \mu, \mu) \leq e^{-2Rt} T_2^2 (f \mu, \mu)
$$
holds for the $T_2$ distance as it does for the $W_2$ distance. The following corollary gives a more precise bound under 
the $CD(R,n)$ condition:

\begin{ecor}
Let $(P_t)_{t \geq 0}$ be our diffusion Markov semigroup on $E$ satisfying a $CD(R,n)$ condition with $R>0$ and $n\geq 1$. Then, in the framework of Corollary \ref{corologsob}, for all $f\in\mathcal F$ and $T\geq0$, 
$$
T_2^2(P_T f\mu,\mu)\leq e^{-2R T}T_2^2(f\mu,\mu) \frac{1}{1 + n \, R \, T_2^2(f\mu,\mu) \, \frac{1-e^{-2RT}}{4 (n-1)^2}}.
$$
\end{ecor}

\begin{eproof}
Taking $g=1$ in \eqref{eq-contraction-T2}, the map 
$
\Lambda(t) = e^{2Rt} T_2^2(P_T f \mu, \mu)$
satisfies
$$
\Lambda'(t) \leq - \frac{nR^2}{2(n-1)^2} e^{-2Rt} \Lambda(t)^2.
$$
Observe indeed that, under the $CD(R,n)$ condition, the measure $\mu$ satisfies the logarithmic Sobolev inequality \eqref{eqLS} with constant $C= \frac{n-1}{2Rn}$ (see \cite{bgl1}), whence a Talagrand inequality \eqref{tala} with constant $4C$ by Corollary \ref{corologsob}. The conclusion follows by time integration.
\end{eproof}

\section{Evolution variational inequalities}\label{sec-EVI}

Evolution variational inequalities (EVI in short) have recently been developed as a connection between curvature conditions $CD(R,\infty)$ (usually in the sense of the commutation of the semigroup and the carr\'e du champ), heat semigroups and the notion of curvature bound introduced by J. Lott, K.-T.~Sturm and C. Villani (see \cite{sturm} and \cite{LV}). We refer to the recent work \cite{AGS13} for a nearly complete picture in this dimensionless setting. However, no dimensional EVI, namely related to a $CD(R,n)$ curvature-dimension condition, were known until M.~Er\nobreak bar, K.~Kuwada and K.-T.~Sturm very recently proved in \cite{EKS13} that $CD(R,n)$ is (roughly) equivalent to
\begin{equation}\label{eks}
\frac d{dt}s_{\frac Rn}\left(\frac 12 W_2(f\mu, H_t g\mu)\right)^2\le -R\,s_{\frac Rn}\left(\frac 12 W_2(f\mu, H_t g\mu)\right)^2+\frac n2\left(1-e^{-\frac1n(\ent{\mu}{f}-\ent{\mu}{H_tg}}\right)
\end{equation}
on the heat semigroup. Here $s_r(x)=\sin(\sqrt{r}x) / \sqrt{r}$ if $r>0$ and $s_r(x)=\sinh(\sqrt{-r}x) / \sqrt{-r}$ if $r<0$. Forgetting for a time the map $s_r$, which is equivalent to $x$ for small $x$, and using that $1-e^{-x}\le x$, this inequality clearly appears to improve the classical EVI obtained under a $CD(R,\infty)$ type condition.

The main goal of this section is twofold. First, we will give a (time integrated) EVI for the Markov transportation distance $T_2$ (rather than the usual Wasserstein distance), under the $CD(R, \infty)$ condition. Then we will see how the possible existence of geodesics can lead to a dimensional EVI, here with a negative corrective term in the spirit of the contraction result in Theorem \ref{thm-contraction-rn}.  As in section \ref{sec-first} we will start with the Euclidean heat equation, i.e. under a $CD(0,n)$ condition, obtaining a dimensional EVI in Wasserstein distance. Assuming the existence of smooth geodesics for $T_2$, we will see that one can obtain a statement under the $CD(R,n)$.

\subsubsection*{Evolution variational inequality for $T_2$ under $CD(R,\infty)$}

\begin{ethm}\label{EVIcdrinfini}
Let $(P_t)_{t \geq 0}$ be our  diffusion Markov semigroup satisfying a $CD(R,\infty)$ condition with $R\in\dR$. Then, for all $f,g\in\mathcal F$,  
$$
T_2^2(f\mu, P_t g\mu)-T_2^2(f\mu,g\mu)\leq -\frac{e^{-2R t}-1+2R t}{2R t}T_2^2(f\mu,g\mu)+2t(\ent{\mu}{f}-\ent{\mu}{P_t g}).
$$
\end{ethm}
 \begin{eproof}
Let $(\rho_s,h_s)$ be an admissible path between $f$ and $g$. Then $(P_{ts}(\rho_s),P_{ts}(h_s-t\rho_s))$ is an admissible path between $f$ and $P_t g$, so
$$
T_2^2(f \mu, P_t g \mu)
\leq
\int_0^1\int \frac{\Gamma(P_{ts}(h_s-t\rho_s))}{P_{ts}(\rho_s)}d\mu ds
$$
\begin{equation}\label{Eq-3termesEVI}
= \int_0^1\! \! \int \frac{\Gamma(P_{ts}(h_s))}{P_{ts}(\rho_s)}d\mu ds
- \, 2t \! \int_0^1 \! \! \int \frac{\Gamma(P_{ts}(h_s-t\rho_s),P_{ts} (\rho_s))}{P_{ts}(\rho_s)}d\mu ds 
- \, t^2 \! \int_0^1\! \! \int \frac{\Gamma(P_{ts}(\rho_s))}{P_{ts}(\rho_s)}d\mu ds.
\end{equation}
By Lemma~\ref{lem4new} under the $CD(R,\infty)$ curvature condition, then for all $s$
$$
\int \frac{\Gamma(P_{ts}(h_s))}{P_{ts}(\rho_s)}d\mu
\leq 
 e^{-2R t s}\int \frac{\Gamma(h_s)}{\rho_s}d\mu
 \leq 
 e^{-2R t s} \big( T_2^2(f\mu,g\mu)+\ep \big)
 $$
 if the admissible path $(\rho_s,h_s)$ is an $\ep$-geodesic, for given $\ep>0$. Moreover
 $$
 \int \frac{\Gamma(P_{ts}(h_s-t\rho_s),P_{ts} (\rho_s))}{P_{ts}(\rho_s)}d\mu
=
\frac{d}{ds}\int P_{ts}(\rho_s)\log P_{ts} (\rho_s)d\mu
$$
by reversibility. Hence, forgetting the last term in \eqref{Eq-3termesEVI} and integrating in $s \in [0,1]$ conclude the argument  by letting $\ep$ go to $0$. 
\end{eproof}

\begin{erem}
Combining the contraction result of Theorem \ref{thm-contraction-rn} under a $CD(R,n)$ condition and this EVI (of course valid under $CD(R,n)$) we may get a contraction type result in the $T_2$ distance for $s<t$, in the spirit of \cite{EKS13}  (see Remark \ref{contts}).
\end{erem}

To obtain dimension dependent bounds under the $CD(R,n)$ condition we will need geodesics. As above for the contraction property, let us see first which additional term coming from the dimension appears for the Euclidean heat equation.

\subsection{A dimensional EVI in Wasserstein distance for the heat equation in $\R^n$}

For simplicity the EVI is here described in its time derivative form, but it may be easily justified by first considering an integrated form of the EVI. 

As in section \ref{sec-first}, the Benamou-Brenier formulation \eqref{eq-BB} is the starting point. Let $(\rho_s,w_s)$ be an admissible path between $f$ and $g$, satisfying the constraint~\eqref{eq-contrainte}. Then $(H_{ts}(\rho_s),R_{ts}(w_s)-t\nabla H_{ts}(\rho_s))$  is an admissible path between $f$ and $H_t g$, and satisfying~\eqref{eq-contrainte}, so by~\eqref{eq-BB},
\begin{equation}
\label{eq-evi-demo}
W_2^2(fdx, H_t gdx)\leq \int_0^1\int\frac{|R_{ts}(w_s)-t\nabla H_{ts}(\rho_s)|^2}{H_{ts}(\rho_s)}dsdx.
\end{equation}

Assume further that $(\rho_s,w_s)$ is a minimizer in the Benamou-Brenier formulation, that is,
$$
W_2^2( fdx,gdx)=\int_0^1\int\frac{|w_s|^2}{\rho_s}dsdx.
$$
The path $(\rho_s)$ is then a geodesic path between $fdx$ and $gdx$ with respect to the Wasserstein distance (see~\cite{dns1} for more details).
In particular inequality~\eqref{eq-evi-demo} is an equality at time $t=0$ and, formally, the time derivative of~\eqref{eq-evi-demo} at $t=0$ implies 
\begin{equation}
\label{eq-avv}
\frac{d}{dt}W_2^2(fdx, H_t gdx)\Big|_{t=0}\leq\frac{d}{dt}\int_0^1\int\frac{|R_{ts}(w_s)-t\nabla H_{ts}(\rho_s)|^2}{H_{ts}(\rho_s)}dsdx\Big|_{t=0}.
\end{equation}
The term on the right-hand side is controlled by the following lemma. 
\begin{elem}
\label{lem2}
Let $(\rho_s,w_s)_{s\in[0,1]}$ be a couple satisfying the constraint~\eqref{eq-contrainte}, where $\rho_s$ is a probability density with respect to Lebesgue measure. Letting 
$$
\Lambda(t)=\int_0^1\int\frac{|R_{ts}(w_s)-t\nabla H_{ts}(\rho_s)|^2}{H_{ts}(\rho_s)}dsdx
$$
for $t \geq 0$, then 
$$
\Lambda'(0)\leq -\frac{2}{n}\int_0^1 s\PAR{\int\frac{w_s\cdot\nabla\rho_s}{\rho_s}dx}^2ds-2\int_0^1\int \frac{w_s\cdot\nabla\rho_s}{\rho_s}dsdx.
$$
\end{elem}
We skip the proof since it is almost the same as  for Lemma~\ref{lem1-lambda}.

\smallskip

Now Lemma~\ref{lem2} and~\eqref{eq-avv} imply
$$
\frac{d}{dt}W_2^2(fdx, H_t gdx) \Big|_{t=0}\leq -\frac{2}{n}\int_0^1 s\PAR{\int\frac{w_s\cdot\nabla\rho_s}{\rho_s}dx}^2ds-2\int_0^1\int \frac{w_s\cdot\nabla\rho_s}{\rho_s}dsdx. 
$$
Letting $\varphi(s)=\int \rho_s\log \rho_sdx$, the relation~\eqref{eq-contrainte} between $w_s$ and $\rho_s$ implies
$$
\varphi'(s)=\int \frac{ w_s\cdot \nabla \rho_s}{\rho_s}dx,
$$
so that
$$
\frac{d}{dt}W_2^2(fdx, H_t gdx)|_{t=0}\leq -\frac{2}{n}\int_0^1 s(\varphi'(s))^2ds-2\int_0^1 \varphi'(s)ds. 
$$
Then the Jensen inequality for the measure $2 \, s \, ds$ and an integration by parts give 
$$
\frac{d}{dt}W_2^2(fdx, H_t gdx)|_{t=0}\leq -\frac{4}{n}\PAR{\varphi(1)-\int_0^1\varphi(s)ds}^2-2( \varphi(1)-\varphi(0)). 
$$
We have obtained the following result :

\begin{eprop}
Let $(H_t)_{t\geq0}$ be the heat semigroup on $\dR^n$. Then, for any probability densities $f$ and $g$ in $\dR^n$ such that $W_2(fdx,gdx)<\infty$,
\begin{equation}
\label{eq-evi-heat}
\frac{1}{2}\frac{d}{dt}W_2^2(fdx, H_t gdx) \Big|_{t=0}\leq -\frac{2}{n} \PAR{\ent{dx}{g}-\int_0^1\ent{dx}{\rho_s}ds}^2+\ent{dx}{f}-\ent{dx}{g}
\end{equation}
where $(\rho_s)_{s\in[0,1]}$ is a geodesic path between $f$ and $g$ for the Wasserstein distance. 
\end{eprop}

At the time being, we have not been able to get a ``geodesic" free version of this dimensional EVI. Note however, once again, that the correction term is quite different in nature from the one obtained for example in \cite{EKS13}. We will see in the next subsection that we obtain results in the same flavor with the $T_2$ distance.

\subsection{A dimension dependent EVI for $T_2$ in the geodesic case}

In this subsection we will assume the existence of smooth geodesics for the $T_2$ distance, 
for any $f,g\in\mathcal F$, there exists  an admissible path  $(\rho_s,h_s)_{s \in [0,1]} \in \mathcal A (f,g)$ such that for all $s\in[0,1]$, 
$$
\int \frac{\Gamma(h_s)}{\rho_s}d\mu= T_2^2(f\mu,g\mu).
$$ 
Therefore, this section is formal, as also we should first consider the integrated form of EVI (which can be obtained, but is however quite difficult to read). We will closely follow the approach used above for the Euclidean heat equation and the Wasserstein distance.

We then consider our  diffusion Markov semigroup $(P_t)_{t \geq 0}$ and, for $f,g \in \mathcal F$, a geodesic path $(\rho_s,h_s)$ between $f$ and $g$ for the associated $T_2$ distance. Then $(P_{ts}(\rho_s)), P_{ts}(h_s-t\rho_s))$ is an admissible path between for $f$ and $ P_t g$. In particular 
\begin{equation}\label{eq-derniere}
T_2^2(f \mu, P_t g\mu)\leq \int_0^1\int \frac{\Gamma(P_{ts}(h_s-t\rho_s))}{P_{ts}(\rho_s)}d\mu ds. 
\end{equation}

We will use the following adaptation of Lemma \ref{lem2} to our setting, which does not need ($\rho_s,h_s)$ to be a geodesic path:
\begin{elem}
\label{lem4}
Let $(\rho_s)_{s \in [0,1]}$ be a smooth path in $\mathcal F$, $(h_s)_{s \in [0,1]}$ be a smooth path in $\mathcal A$, and
$$
\Lambda(t)=\int_0^1\int\frac{\Gamma(P_{ts}(h_s-t\rho_s))}{P_{ts}(\rho_s)}d\mu ds 
$$
for $t \geq 0$. Then, under the $CD(R,n)$ condition, 
$$
\Lambda'(0)\leq -2 R \int_0^1s\int\frac{\Gamma(h_s)}{\rho_s}d\mu ds-\frac{2}{n}\int_0^1s\PAR{\int \frac{\Gamma(h_s,\rho_s)}{\rho_s}d\mu}^2ds-2\int_0^1\int \frac{\Gamma(h_s,\rho_s)}{\rho_s}d\mu ds.
$$
\end{elem}

\bigskip

Since $(\rho_s, h_s)$ is a geodesic path, then \eqref{eq-derniere} is an equality at time $t=0$, so, taking the time derivative at $t=0$ and using Lemma~\ref{lem4}, 
\begin{eqnarray*}
\frac{1}{2} \frac{d}{dt}T_2^2(f\mu, P_t g\mu) \Big|_{t=0}&\leq&-R \int_0^1s\int\frac{\Gamma(h_s)}{\rho_s}d\mu ds\\&&\qquad-\frac{1}{n}\int_0^1s\PAR{\int \frac{\Gamma(h_s,\rho_s)}{\rho_s}d\mu}^2ds-\int_0^1\int \frac{\Gamma(h_s,\rho_s)}{\rho_s}d\mu ds.
\end{eqnarray*}
Again $(\rho_s, h_s)$ is a geodesic path, so
$$
\int\frac{\Gamma(h_s)}{\rho_s}d\mu =T_2^2(f\mu, g\mu) 
$$
for all $s$. Letting again $\varphi(s)=\int \rho_s\log\rho_sd\mu$, the inequality may then be rewritten as 
$$
\frac{1}{2} \frac{d}{dt}T_2^2(f\mu, P_t g\mu) \Big|_{t=0}\leq-\frac{R}{2} \, T_2^2(f\mu,g\mu)-\frac{1}{n}\int_0^1s \, \varphi'(s)^2 \, ds-\int_0^1\varphi'(s) ds.
$$
Hence, by the Jensen inequality for the measure $2sds$ and an integration by parts, the $CD(R,n)$ condition and the existence of geodesics ensure  the following dimensional EVI:
\begin{equation}\label{dimEVI}
\frac{1}{2}\frac{d}{dt}T_2^2(f\mu, P_t g\mu) \Big|_{t=0}\leq -\frac{R}{2}T_2^2(f\mu,g\mu)\\-\frac{2}{n}\PAR{\ent{\mu}{g}-\int_0^1\ent{\mu}{\rho_s}ds}^2+\ent{\mu}{f}-\ent{\mu}{g},
\end{equation}
where $(\rho_s)$ is a geodesic path between $f\mu$ and $g\mu$ for the $T_2$ distance.

\section{$\Phi$-entropies versus usual entropy}
\label{sec-gene}

There are many ways of extending the Markov transportation distance. Here we present the one associated with $\Phi$-entropies, well adapted to the $\Gamma_2$-calculus.  For the Wasserstein distance this generalization has been formulated in~\cite{dns1,dns2}. 

Let again $(P_t)_{t \geq 0}$ be a diffusion Markov semigroup with invariant measure $\mu$. Let $\xi$ be a $\mathcal C^2$ positive function on $(0,+\infty)$ with $1/\xi$ concave. Let also  
$$
\entf{\mu}^\Phi(f)=\int\Phi(f)d\mu-\Phi \Big(\int fd\mu \Big)
$$
 be the $\Phi$-entropy of a positive map $f$, with $\Phi''=\xi$. The $\Phi$-entropies have been studied for instance in~\cite{bakrysaintflour,chafai,arnolddolbeault,bolleygentil}. By analogy with Definition \ref{defiT2}, for $\Phi(x) = x \, \log x$ and $\xi (x) = 1/x$:

\begin{edefi}
For $f,g\in\mathcal F$ we let
$$
T_\xi(f\mu,g\mu)=  \inf \Big( \int_0^1\int {\Gamma(h_s)}{\xi(\rho_s)}d\mu ds \Big)^{1/2} 
$$ 
where the infimum runs over all admissible paths $(\rho_s,h_s)_{s \in [0,1]} \in\mathcal A(f,g)$.
\end{edefi}
\bigskip

For a general map $\xi$ this distance shares the same properties of existence of $\varepsilon$-geodesics and tensorization as the distance $T_2$, which can be proved as in section~\ref{sec-defusu}. For instance:

\begin{eprop}[Tensorization]
Let $\xi$ be a $\mathcal C^2$ positive function on $(0, + \infty)$ with $1/\xi$ concave. Let $(P_t^i)_{t\geq0}$, $i\in\{1,\cdots, N\}$ be $N$ Markov semigroups on  compact connected Riemannian manifolds $E_i$ with probability measure  $\mu_i$, with generators $L_i$ and carr\'es du champ $\Gamma_i$ as in Section~\ref {sec-defdef}.  Let  $P_t=\otimes_{i=1}^N P_t^i$ on the product space $E=\times_{i=1}^N E_i$, as in Proposition~\ref{prop-tensorisation}.

 Then, for any densities $f(x)=\prod_{i=1}^N f_i(x_i)$ and $g(x)=\prod_{i=1}^N g_i(x_i)$ ($x=(x_1,\cdots,x_N)$) in $\mathcal F$, 
\begin{equation}
\label{eq-tensoxi}
T_\xi^2(f\mu,g\mu)\geq \sum_{i=1}^NT_{\xi}^2(f_i\mu_i,g_i\mu_i).
\end{equation}
\end{eprop}

\begin{eproof}
The argument follows the proof of  Proposition~\ref{prop-tensorisation}. The key bound \eqref{eq-jensen} is replaced by
$$
\int \Gamma_1(h_s)\xi(\rho_s)d\mu_2\geq \Gamma_1(h_s^1)\xi(\rho_s^1).
$$
This is a consequence of the definition~\eqref{eq-gamma} of the carr\'e du champ and of the fact that the map $(x,y)\mapsto x^2\xi(y)$ is convex under our assumption on $\xi$, see~\cite{chafai}.
\end{eproof}

\subsubsection*{Contraction and evolution variational inequalities}

\begin{ethm}
Let $(P_t)_{t \geq 0}$ be our diffusion Markov semigroup satisfying  a $CD(R,\infty)$ condition with $R\in\dR$. Then for any $f,g\in\mathcal F$ and $t\geq0$, the contraction property 
\begin{equation}
\label{eq-phi-contraction}
T_\xi^2(P_t f\mu,P_t g\mu)\leq e^{-2R t}T_\xi^2(f\mu,g\mu)
\end{equation}
holds, as well as  the Evolution Variational Inequality
\begin{equation}
\label{eq-phi-evi}
T_\xi^2(f\mu, P_t g\mu)-T_\xi^2(f\mu,g\mu)\leq -\frac{e^{-2R t}-1+2R t}{2R t} \, T_\xi^2(f\mu,g\mu)+2t \, (\entf{\mu}^\Phi(f)-\entf{\mu}^\Phi(P_t g)),
\end{equation}
where $\Phi''=\xi$. 
\end{ethm}

\begin{eproof}
The proof follows the idea of the classical case of Theorems \ref{thm-contraction-rn} and \ref{EVIcdrinfini}. It uses the relation
$$
\int_0^1\int \Gamma(P_{ts}(h_s-t\rho_s),P_{ts} (\rho_s))\xi(P_{ts}(\rho_s)) \, d\mu \, ds=\entf{\mu}^\Phi(P_t g)-\entf{\mu}^\Phi(f).
$$
for any admissible path $(\rho_s,h_s)$ between $f$ and $g$, and in particular
$$
\int_0^1\int \Gamma(\rho_s,h_s) \xi (\rho_s) \, d\mu \, ds=\entf{\mu}^\Phi(g)-\entf{\mu}^\Phi(f),
$$
and the following lemma.
\end{eproof}

\begin{elem}
\label{lem-phifonda}
Let $f,g\in\mathcal A $ with $g>0$ and let  
 $$
 \Lambda(t)=\int {\Gamma(P_tf)}{\xi(P_t g)}d\mu
 $$
for $t \geq 0$. Then, under the curvature condition $CD(R,\infty)$, 
 \begin{equation}
 \label{eq-phi-derive}
\Lambda'(t)\leq -2 R\Lambda(t)-\int \xi^2(P_tg)\PAR{-\frac{1}{\xi}}''\!\!\!(P_tg)\Gamma(P_tf)\Gamma(P_tg)d\mu
 \end{equation}
\end{elem}

\begin{eproof}
We only briefly check the proof since it follows the one of Lemma~\ref{lem4new}. For any $t\geq0$, 
$$
\Lambda'(t)=\int\Big[2\Gamma(P_tf,LP_tf)\xi(P_tg)+\Gamma(P_t f)\xi'(P_tg)LP_tg\Big]d\mu.
$$
In the notation $G=P_tg$ and $F=P_tf$, the invariance property $\int L\big[\Gamma(F)\xi(G)\big]d\mu=0$ and the diffusion property of $L$ give that 
$$
\int\Gamma(F)\xi'(G)LGd\mu=-\int\Big[2\xi'(G)\Gamma(\Gamma(F),G)+\xi''(G)\Gamma(G)\Gamma(F)+\xi(G)L\Gamma(F)\Big]d\mu.
$$
Hence, using the definition of the $\Gamma_2$ operator,  
\begin{multline*}
\Lambda'(t)=-2\int\xi(G)\Big[\Gamma_2(F)+\Gamma(\Gamma(F),\log \xi(G))+\frac{\xi''(G)\xi(G)}{2\xi'(G)^2}\Gamma(F)\Gamma(\log \xi(G))\Big]d\mu\\
=-2\int\xi(G)\Big[\Gamma_2(F)+\Gamma(\Gamma(F),\log \xi(G))+\Gamma(F)\Gamma(\log \xi(G))\Big]d\mu
-\int \xi^2(G)\PAR{-\frac{1}{\xi}}''\!\!\!(G)\Gamma(F)\Gamma(G)d\mu.
\end{multline*}
Then Lemma~\ref{lem-technique} for $n = \infty$, applied to $f=F$ and $g=\log\xi(G)$, implies inequality~\eqref{eq-phi-derive}. 
\end{eproof}

\subsubsection*{The particular case of power functions}
~

Poincar\'e and logarithmic Sobolev inequalities belong to the family of $\Phi$-entropy inequalities, namely for $\Phi(x) = x^2/2$ and $\Phi (x) = x \, \log x$ respectively (see \cite{chafai}). An interpolation family of inequalities between them consist in the Beckner inequalities, for $\Phi_p(x)=\frac{x^p}{p(p-1)}$. It has been proved in~\cite{arnolddolbeault,bolleygentil} how to refine these Beckner inequalities under the curvature-dimension condition $CD(R,\infty)$. In the same way, the contraction inequalities proved in~\cite{dns2} and in \eqref{eq-phi-contraction} for a general $\Phi$ can be made more precise for these power functions, as follows.

For $p\in(1,2)$ we let $\xi_p(x)=x^{p-2}$ and $\Phi_p(x)=\frac{x^p}{p(p-1)}$ for $x>0$, so that $\Phi_p'' = \xi_p$.

\begin{ethm}[Refined contraction inequality]
Let $(P_t)_{t \geq 0}$ be our  diffusion Markov semigroup satisfying  a $CD(R,\infty)$ condition with $R\in\dR$. Then, for any $f,g\in\mathcal F$ and $t\geq0$,
\begin{multline}
\label{eq-derdesder}
T_{\xi_p}^2(P_tf\mu, P_tg\mu)\leq 
e^{-2Rt}T_{\xi_p}^2(f\mu,g\mu)\\
-4 \, \frac{2-p}{p^2(p-1)}\int_0^te^{-2R(t-u)}\PAR{\sqrt{\int(P_uf)^pd\mu}-\sqrt{\int(P_ug)^pd\mu}}^2du.
\end{multline}
\end{ethm}
\begin{eproof}
Let $(\rho_s, h_s)$ be an admissible path between $f$ and $g$, and $\Lambda(t,s) = \int \Gamma(P_t(h_s)) \xi_p(P_t (\rho_s)) d\mu$. Then inequality~\eqref{eq-phi-derive} for $h_s$ and $\rho_s$ writes 
$$
\partial_t\Lambda(t,s)\leq-2R\Lambda(t,s)-(2-p)(p-1)\int (P_t(\rho_s))^{p-4}\Gamma(P_t(\rho_s))\Gamma(P_t(h_s))d\mu. 
$$
But
$$
\int (P_t(\rho_s))^{p-4}\Gamma(P_t(\rho_s))\Gamma(P_t(h_s))d\mu\geq \frac{\Big(\int (P_t(\rho_s))^{p-2}\Gamma(P_t(\rho_s),P_t(h_s))d\mu\Big)^2}{\int P_t(\rho_s)^pd\mu}
$$
by the Cauchy-Schwarz inequality,  so
$$
\partial_t\Lambda(t,s)\leq-2R\Lambda(t,s)-\frac{2-p}{p}\frac{[\partial_s\phi(t,s)]^2}{\phi(t,s)}
$$
where 
$$
\phi(t,s)= \frac{1}{p(p-1)}\int (P_t(\rho_s))^pd\mu.
$$
Integrating over $s\in[0,1]$ and applying the Gronwall inequality in $t$, we obtain 
$$
\int_0^1\Lambda(t,s)ds\leq e^{-2Rt}\int_0^1\Lambda(0,s)ds-\frac{2-p}{p}\int_0^t\int_0^1e^{-2R(t-u)}\frac{[\partial_s\phi(u,s)]^2}{\phi(u,s)}dsdu. 
$$
But
$$
\int_0^1\frac{[\partial_s\phi(u,s)]^2}{\phi(u,s)}ds\geq \PAR{\int_0^1\frac{\partial_s\phi(u,s)}{\sqrt{\phi(u,s)}}ds}^2=4(\sqrt{\phi(u,1)}-\sqrt{\phi(u,0)})^2
$$
again by the Cauchy-Schwarz inequality. The result follows by optimizing over $(\rho_s, h_s)$.
\end{eproof}

In the limit case where $p=2$, the improvement in the contraction inequality disappears in~\eqref{eq-derdesder}, as observed in the refined Beckner inequalities of~\cite{arnolddolbeault,bolleygentil}. Morever this improvement goes to $0$ when $p$ goes to 1, hence recovering the classical contraction inequality \eqref{eq-contraction-T2}-\eqref{eq-phi-contraction} under the  curvature condition $CD(R,\infty)$. 



\bigskip

\newcommand{\etalchar}[1]{$^{#1}$}
{\footnotesize{

}

\noindent{\bf Acknowledgements.} The authors warmly thank Guillaume Carlier, Nicola Gigli, Christian L\'eonard and Karl-Theodor Sturm for enlightening discussion. They are grateful to the referee for a careful reading of the manuscript and most relevant
comments and questions which helped improve the presentation of the paper. This research was supported  by the French ANR STAB project.

\end{document}